\documentclass[12pt,a4paper]{article}

\usepackage{cite}
\usepackage[utf8]{inputenc}
\usepackage[english]{babel}
\usepackage{amsmath}
\usepackage{mathrsfs}
\usepackage{amsfonts}
\usepackage{amsthm}
\usepackage{amssymb}
\usepackage{makeidx}
\usepackage{graphics}
\usepackage{float}
\usepackage{subfloat}
\usepackage{subfigure}
\usepackage{wrapfig}
\usepackage[usenames, dvipsnames]{color}
\usepackage{tcolorbox}
\usepackage{hyperref}
\definecolor{link}{rgb}{0.1,0.1,0.9}
\hypersetup{colorlinks=true,linkcolor=link,citecolor=link,urlcolor=link,linktocpage}

\usepackage[left=1.8cm,right=1.8cm,top=2cm,bottom=2cm]{geometry}
\usepackage[capitalize]{cleveref}
\usepackage{tikz}
\usetikzlibrary{shapes}
\usetikzlibrary{plotmarks}
\usepackage[normalem]{ulem}

\makeatletter
\def\Let@{\def\\{\notag\math@cr}}
\makeatother

\usepackage{authblk}

\newcommand{\be}{\begin{equation}}
\newcommand{\ee}{\end{equation}}
\newcommand{\ben}{\begin{eqnarray}\displaystyle}
\newcommand{\een}{\end{eqnarray}}

\newcommand{\bea}[1]{\begin{eqnarray}\label{#1} }
\newcommand{\eea}{\end{eqnarray}}

\newcommand{\half}{\frac{1}{2}}

\newtheorem{Theorem}{Theorem}[section]
\newtheorem{lemma}{Lemma}[section]
\newtheorem{prop}{Proposition}[section]
\newtheorem{conj}{Conjecture}[section]

\def\boxempty{{\,\lower0.9pt\vbox{\hrule \hbox{\vrule height 0.25 cm
\hskip 0.25 cm \vrule height 0.25 cm}\hrule}\,}}

\begin{document}
\begin{titlepage}

\title{
{\Huge\bf Shifted second moment of the Riemann}\\
{\Huge\bf zeta function and a Fourier type kernel }
}

\bigskip\bigskip\bigskip\bigskip\bigskip

\author{{\bf Parikshit Dutta}${}^{1}$\thanks{{\tt parryparikshit@gmail.com}},
                    {\bf Debashis Ghoshal}${}^2$\thanks{{\tt dghoshal@mail.jnu.ac.in}}, \\  
                        {\bf and Krishnan Rajkumar}$^3$\thanks{{\tt krishnan.rjkmr@gmail.com}}\\
\hfill\\        
${}^1${\it Asutosh College, 92 Shyama Prasad Mukherjee Road,}\\
{\it Kolkata 700026, India}\\      

\hfill\\
${}^2${\it School of Physical Sciences, Jawaharlal Nehru University,}\\
{\it New Delhi 110067, India}\\

\hfill\\
${}^3${\it School of Engineering, Jawaharlal Nehru University,}\\
{\it New Delhi 110067, India}
}

\bigskip\bigskip\bigskip\bigskip

\date{%
%
\bigskip\bigskip
\begin{quote}
\centerline{{\bf Abstract}}
{\small
We compute the second moment of the Riemann zeta function for shifted arguments over a domain that extends the ones in the literature. 
We use the Riemann-Siegel formula for the error term in the approximate functional equation and take the products of all the terms into 
account. We also show that, as a function of imaginary shifts on the critical line, the the second moment behaves like a Fourier-Cauchy 
type kernel on a class of functions. This is reminiscent of orthogonal functions.}
\end{quote}
}

\bigskip


\end{titlepage}
\thispagestyle{empty}\maketitle\vfill \eject

\tableofcontents

\section{Introduction}\label{sec:Introd}
The study of the moments of the Riemann zeta function $\zeta(s)$ in the critical strip $0\leq \Re s \leq 1$, with its long history dating back 
to the work of Hardy-Littlewood, remains an active area of research. For an elucidating account of this field and its connections with the 
value distribution of $\zeta(s)$, see the reviews by Soundararajan \cite{Sound:2008Moments,Sound:2021Distribution} (see also 
Refs.\cite{titchmarsh1986theory,ivic2013riemann,iwaniec2004analytic}).

Let $I_k(\sigma, T)$ denote the $2k$-th moment of the Riemann zeta function on the half plane $\sigma\geq\half$ be defined as
\begin{align*}
	I_k(\sigma, T) = \int_{T_0}^T \left|\zeta(\sigma+ i t)\right|^{2k} dt
\end{align*}
where $T> T_0 >0$ and $T_0$  any constant. The following result may be extracted from 
Refs.\cite{titchmarsh1986theory} (Ch.~VII) and \cite{ivic2013riemann} (Ch.~1). 
\begin{Theorem}\label{ivic}
For $0 < k \leq 2$ and $\sigma>\half$, it is known that
\begin{align}\label{iviceq}
	I_k(\sigma,T) = \left(\sum_{n=1}^{\infty} \frac{d_k^2(n)}{n^{2\sigma}}\right) T + o(T)
\end{align}
For $k=1$, the error term can be improved to $ O(T^{2-2\sigma+\epsilon})$ in the region $\half<\sigma<1$. 
For any $k>0$, Eq.\eqref{iviceq} holds with error term $O(T^{2-\sigma+\epsilon})$ in $\sigma>1$. 
\end{Theorem}
In particular, on the critical line $\sigma = \tfrac{1}{2}$, several authors have derived conditional and unconditional bounds for 
$I_k(\tfrac{1}{2},T)$ (see Refs.\cite{Sound:2008Moments,Sound:2021Distribution} for a comprehensive set of references). 
The present state of the art may be summarised as follows.
\begin{Theorem} 
Let $k > 0$ and $T \geq e$ be real numbers. There are positive constants $a_k$ and $b_k$ such that
\begin{align}\label{lubounds}
	a_k T(\log T)^{k^2} \leq  I_k(\tfrac{1}{2},T) \leq b_k T(\log T)^{k^2}
\end{align}
in which the lower bound holds unconditionally for all $k$, while the upper bound holds unconditionally for $k \leq 2$ and, conditional on 
the Riemann hypothesis, for all $k > 2$.
\end{Theorem}

{}From the bounds in Eq.\eqref{lubounds}, it is reasonable to surmise that $I_k(\tfrac{1}{2},T) \sim c_k T \log^{k^2} T$ (where $c_k$ is some 
constant) for all $k > 0$. However, this is known only for $k=1$ \cite{HardyLittle:AM1918} and $k=2$ \cite{Ingham:MV1928}. For any $k>0$, 
a conjecture due to Keating and Snaith \cite{KeatingSnaith:2000RMTZeta}, states that
\begin{conj}\label{KeatingSnaith}
	\begin{align*}
		I_k(\tfrac{1}{2},T) \sim \dfrac{\Gamma(k^2+1)G(k+1)^2}{G(2k+1)} \sum_{n\leq T} \frac{d_k^2(n)}{n}
	\end{align*}
where $G(n)=\prod_{i=0}^{n-2}i!$ for $n\in \mathbb{N}$, is the Barnes $G$-function.	
\end{conj}
In the above, the notation used is consistent with that in Theorem \ref{ivic}, as opposed to the original one in Ref.\cite{KeatingSnaith:2000RMTZeta} 
where the sum on the right is 
replaced by $a(k) (\log T)^{k^2}$. The conjecture is based on an analogy between the average value of $|\zeta(\tfrac{1}{2} + i t)|^{2k}$ in the 
range $[T_0,T]$ and the average of the $2k$-th power of the absolute value of the characteristic polynomial of a random ensemble of unitary 
$N\times N$ matrices, with $N  \sim \log \frac{T}{2 \pi}$. This connection goes back to the work of 
Montgomery \cite{Montgomery:1972, Anecdote:2003}. Analogous conjectures exist for Dirichlet $L$-functions \cite{KeatingSnaith:2000RMTL}.

The shifted integral moments 
\begin{align*}
	I(a_1,\cdots,a_k,a_{k+1},\cdots,a_{2k}; T) = \int_{T_0}^T dt\, \prod_1^k \zeta(\tfrac{1}{2}+a_i+ i t) \zeta(\tfrac{1}{2}-a_{k+i}- i t)
\end{align*}
are also of interest. For $k=1$ and complex values of $a_1 = a$ and $a_2 = b$, Ingham \cite{Ingham:MV1928} gave the complete 
asymptotic expression as
\begin{align}\label{ingham}
	I (a,b; T) = \int_{T_0}^T d \tau\, \left( \zeta(1+a-b) + \psi(a,b;\tau) \zeta(1+b-a) \right) + O\left(T^{\half + \frac{b_R-a_R}{2}} \log^2 T\right)
\end{align}
where $\psi(a,b; t) = \left(\tfrac{t}{2\pi}\right)^{b-a}$, $T_0>0$, $a_R=\Re a$, $a_I=\Im a$, $b_R=\Re b$, $b_I=\Im b$ and the shifts $a,b$ lie 
in the rectangular region $\left\{ z : |\Re z| \leq \tfrac{1}{2},\; |\Im z|< T_0\right\}$. This result was generalized by Bettin \cite{Bettin:2011Second} 
to show that Eq.\eqref{ingham} holds with the function $\psi(a,b;\tau) = \chi(\tfrac{1}{2} + a + i \tau) \chi(\tfrac{1}{2} - b - i \tau)$ for $a$ and $b$ 
in the region $\left\{ z : |\Re z| \ll \frac{1}{\log T},\; |\Im z|\ll T\right\}$ and $\chi(s)$ as in Eq.\eqref{eq:DefChi}. In other words, this allowed 
Ingham's result to hold for unbounded shifts along the imaginary axis, although the real parts of the shifts are asymptotically close to zero.

Generalizing the above for $k\geq 2$, Conrey {\em et al} \cite{conrey2005integral} proposed the following conjecture.
\begin{conj}\label{conrey} 
For an even number of complex parameters $(a_i, a_{k+i})$, $i= 1, 2, \cdots, k$, such that $-\half < \Re{a_i}, \Re{a_{k+i}} <\half$ for all 
$i$, the $2k$-th moment of the zeta function with shifted arguments is
\begin{equation*}
		I \left(\left\{a_i, a_{k+i}\right\}; T\right) = \int_{T_0}^T dt \, W_k \left( \left\{a_i\right\}; t\right) \left(1 + O(t^{-\half+\epsilon})\right)
\end{equation*}
where 
\begin{align}\label{W}
	W_k  \left(\left\{a_i, a_{k+i}\right\}; t \right)  &= 
	\sum_{\sigma\in G}  \left(\frac{t}{2\pi}\right)^{\half \sum_{i} \left((a_{\sigma(i)} - a_i) - (a_{\sigma(k+i)}  - a_{k+i})\right)}
		 A_k \left(\left\{a_{\sigma(i)}, a_{\sigma(k+i)}\right\}\right) \\
        &\qquad\times\,\prod_{i,j=1}^k \zeta(1 + a_{\sigma(i)} - a_{\sigma(k+j)})
\end{align}
where $G = \left\{\sigma \in S_{2k} :  \sigma(1) < \cdots < \sigma(k)  \text{ and }  \sigma(k+1) < \cdots < \sigma(2k)\right\}$ is a subgroup 
of the permutation group $S_{2k}$ and $A_k$ is an Euler product involving the parameters.
\end{conj}
The approach in Ref. \cite{conrey2005integral} is based on a heuristic `recipe' considering only the diagonal terms in the  products of 
Dirichlet series appearing in the approximate functional equation for the zeta function (see \cref{sec:FuncEqn}) and {\em ignoring} the 
off-diagonal contributions and the error terms in all the resulting products. They also compared this with the analogous calculation from 
random matrix theory (RMT), which gave an analogous sum with the following substitutions: the argument $\frac{t}{2\pi}$ is replaced by 
$e^{N}$ and the $\zeta(1+z)$ terms are replaced by $(1-e^z)^{-1}$. However, the factors $A_k$, which contain arithmetic information are 
trivial on the RMT side. This means that the calculations would match (at the leading order) only at the poles, when $a_1 = \cdots = 
a_{2k} = 0$, and take us back to Conjecture \ref{KeatingSnaith}.

It may be noted that for $k=1$, Conjecture \ref{conrey} agrees with Ingham's result Eq.\eqref{ingham} up to the leading terms, since 
$A_2(a,b)=1$ in that case. However, the conjectured $O(T^{\half+\epsilon})$ error term seems rather tight in comparison with 
Eq.\eqref{ingham}, which might be expected given the known error terms in the $k=1$ case in Theorem \ref{ivic}.

\bigskip

In this paper, we present two main results in Secs. \ref{sec:CompMom2} and \ref{sec:ZKernel}, respectively. The first is a generalization 
of the result of Bettin \cite{Bettin:2011Second} that Eq.\eqref{ingham} holds in a rectangular strip with asymptotically large vertical shifts.
\begin{Theorem}\label{thm1} 
For $a,b$ in the region $\left\{z : |\Re z| < \epsilon,\, |\Im z| < T\right\}$, the following holds uniformly.
	\begin{align*}
		\int_{-T}^T dt \ \zeta(\tfrac{1}{2}+a+ i t)\zeta(\tfrac{1}{2}-b- i t)
	 &= \int_{-T}^T dt\, \left( \zeta(1+a-b)+\psi(a,b,t)\zeta(1+b-a) \right) \notag\\
	 &\qquad +O(T^{\half}|a_I-b_I|^{\half+2\epsilon})+O(T^{\frac{1}{2}+3\epsilon})
	\end{align*}
where $\psi(a,b,t) = \chi(\tfrac{1}{2}+a+i t) \chi(\tfrac{1}{2}-b-i t)$.
\end{Theorem}
Our approach is to start with the approximate functional equation with a precise form for the error term (using the Riemann-Siegel formula) 
and analyse the products of all the terms from the approximate functional equation, {\em including} the off-diagonal contributions as well as 
the integrals involving the error terms. This approach improves upon Ref.\cite{Ingham:MV1928}, since we use the Riemann-Siegel formula 
(which was not available to Ingham) to get much stronger estimates when $a,b$ are {\em not} confined to a bounded rectangle. The reliability 
of the proof by Bettin \cite{Bettin:2011Second} seems to be in question since the main terms are obtained by the use of Eqs. (7.4.2) and (7.4.3) 
(in p.8) of Ref. \cite{titchmarsh1986theory}. This works only when the integral in $I(a,b,T)$ is from $-T$ to $T$. The error term (involving
$a_I - b_I$ in particular) may also need a reassessment.
 
\medskip

The second result is an application of Theorem \ref{thm1}.  Recall that the Fourier-Cauchy kernel 
\begin{equation*}
	\frac{1}{2T}\int_{-T}^T dt\, e^{it(\alpha - \beta)} 
	= \frac{\sin T(\alpha - \beta)}{T (\alpha -\beta)}
\end{equation*}
is the Dirac delta function (in the distributional sense) in the limit $T\to\infty$. Let us define, in an analogous fashion, the kernel
\begin{equation}
	\mathcal{K}_\zeta(\alpha,\beta;T) = \frac{1}{2T} \int_{-T}^T dt\ \zeta(\tfrac{1}{2}+i \alpha + i t)\zeta(\tfrac{1}{2}-i \beta - i t)
	\label{eq:Zeta2MomKernel}
\end{equation}
\begin{Theorem}\label{thm2} 
For a class of functions $f(\alpha)$ that are meromorphic in $\mathbb{C}$, holomorphic in $\Im \alpha\geq 0$, $f(\alpha), f(\alpha)|
\alpha|^{\frac{1}{2}+2\epsilon} \in L^1(\mathbb{R})$ and satisfy the asymptotic conditions $f(z) = o(|z|^{-|z|})$ in the UHP and 
$f(z)=o(e^{-c|z|})$ in the LHP as $|z|\rightarrow\infty$, the following holds.
	\begin{equation}
		\lim_{T\to\infty} \int_{-T}^T d\alpha\, f(\alpha)\, \mathcal{K}_\zeta(\alpha,\beta;T) = 2\pi f(\beta)
		\label{eq:Zeta2KernTestFn}
	\end{equation}
\end{Theorem}
\noindent
This property of the moment was observed in Refs.\cite{Atkinson:1948MVRZ,Kosters:2008RZSine}, however, the range of the integrals 
there are as in Eq.\eqref{ingham}.

\medskip

To complete the analogy, we consider the following kernel in RMT (see also Ref.\cite{Kosters:2008RZSine})
\begin{equation}
	\mathcal{K}_{\text{CUE}}(x,y; N) = \int_{U(N)}dM \:\Lambda_{M}(e^{-ix})\Lambda_{M^{\dagger}}(e^{iy})
	\label{eq:CUEKernel}
\end{equation}
\begin{Theorem}\label{thm3} 
For a class of functions $f(\alpha)$ that are meromorphic in $\mathbb{C}$, holomorphic in the open unit disc $|z| < 1$ and 
satisfy the asymptotic condition $f(z) = o(|z|^{N+1})$ as $|z|\rightarrow\infty$, the following holds.
	\begin{equation}
		\lim_{N\to\infty}\int_{0}^{2\pi}dx\,f(e^{ix})\, \mathcal{K}_{\text{CUE}}(x,y; N) = 2 \pi f(e^{iy}) \label{eq:CUEtheorem}
	\end{equation}
\end{Theorem}

In this context, a result analogous to Theorem \ref{thm3}, involving the Bessel functions, may be noted.
\begin{Theorem} (Hankel \cite{watson1995treatise})\label{thm4} 
For any (real valued) function $F(r)$ on $\mathbb{R}^+$ that satisfies the boundedness condition 
$\displaystyle{\int_{0}^{\infty} |F(r)| \sqrt{r}\, dr} < \infty$,
\begin{equation}\label{HankelT}
 	\int_{0}^{\infty} du\, u \int_{0}^{\infty} dr\, r F(r) J_{\nu}(ur) J_{\nu}(ut) = \frac{1}{2} \left(F(t+) + F(t-)\right),
	\quad \text{for } \nu>-\frac{1}{2}
\end{equation}
\end{Theorem}
For a detailed review and proof, see Ref.\cite{watson1995treatise}. We would like to point out the analogy between Hankel transform and 
Theorems \ref{thm2} and \ref{thm3}. In particular, this can be interpreted as an orthogonality relation on the space of test functions that 
satisfy appropriate conditions. This is akin to several other transforms (e.g., Mehler-Fock, Kontorovich-Lebedev and cosine), which employ 
orthogonal functions such as the Legendre function $P_{it-\frac{1}{2}}(x)$ or the modified Bessel functions of second kind $K_{it}(x)$. Thus, 
Theorem \ref{thm2} is suggestive of the zeta function on the critical line $\zeta\left(\half + it\right)$ being similar to an orthogonal function. 

An important property for these functions, which are in close resemblance to orthogonal {\em polynomials}, is the reality and simplicity 
of the zeroes as a function of the index. While this property of orthogonal polynomials is easy to establish, the same is not the case for 
orthogonal functions. For instance, the reality and simplicity of the zeroes of the modified Bessel functions $K_{it}(x)$ (with respect to the 
index $t$) was proved by Polya \cite{Polya} (see also \cite{Hejhal:1994}). This suggests that a general proof of the properties of the zeroes 
of orthogonal functions may also be instructive for the zeta function.

\medskip

Finally, let us mention that the approach adopted in this paper very likely applies {\em mutatis mutandis} to the Dirichlet $L$-functions.

\section{Approximate functional equation and other results}\label{sec:FuncEqn}
In this section, in preparation of our main theorems, we recall and adapt a few important results from the literature.

In Ref.\cite{titchmarsh1986theory} (Theorem 4.16, with $N=3$) Titchmarsh proves that for $0 \le \sigma \le 1$, $t>C>0$ and 
$q=\lfloor t/2\pi\rfloor$, the zeta function satisfies the approximate functional equation
\begin{align}\label{afe}
\zeta(s) = \sum_{n\le q} n^{-s} + \chi(s) \sum_{n\le q} n^{s-1} + E(s)
\end{align}
where 
\begin{align}
\chi(s) = 2 (2\pi)^{s-1} \Gamma(1-s) \cos(\tfrac{\pi}{2}(1-s)) = - 2(2\pi)^{s-1}s\Gamma(-s)\sin \frac{\pi s}{2}
\label{eq:DefChi}
\end{align}
and the error term $E(s)$ is the truncation of an asymptotic series in the Riemann-Siegel formula which simplifies to
\begin{align}\label{siegel-term}
	&E(s)=(-1)^{q-1} e^{- \frac{i\pi}{2}(s-1) - \frac{i}{2} t - \frac{i\pi}{8}} (2\pi t)^{\frac{1}{2}(s-1)} \Gamma(1-s) 
	\phi(t) + O(t^{-\frac{3}{4}}) 
\end{align}
where
\begin{align*}
	\phi(t) = \frac{\cos\left(t - (2q+1)\sqrt{2\pi t} - \frac{\pi}{8}\right)}{\cos \sqrt{2\pi t}}
\end{align*} 
We need a slightly stronger version of a result of Titchmarsh \cite{titchmarsh1986theory} (Lemma, p.143), the proof of which remains 
the same.
\begin{lemma}\label{titch_lemma} (Titchmarsh)  
 \begin{align*}
 	 \displaystyle{\int_{\gamma-iT}^{\gamma+iT}\!\! \frac{ds}{2\pi i}\, \chi(1-s)x^{-s}} = 
	 \begin{cases}
	 2 \cos( 2\pi x ) + O\left(\frac{T^{\gamma-\half}}{x^{\gamma} \log \frac{T}{2 \pi x}} \right) + 
	 O\left(\frac{T^{\gamma-\half} \log T}{x^{\gamma}}\right) &\text{if } x < \frac{T}{2 \pi} \text{ and } 0 < \gamma < 1\\
	 O\left(\frac{T^{\gamma - \half}}{x^{\gamma} \log \frac{2 \pi x}{T}}\right) + O\left(\frac{T^{\gamma - \half}}{x^{\gamma}}\right)
	 &\text{if } x > \frac{T}{2 \pi} \text{ and } \gamma > \frac{1}{2}
	 \end{cases}
\end{align*} 
 \end{lemma}
\begin{prop}\label{hilbert}
	For $N \in \mathbb{N}$ and $a_i,b_i \in \mathbb{C}$, $i=1, 2, \cdots, N$, the following versions of Hilbert's inequality hold (see 
	Ref.\cite{ivic2013riemann}, Sec 5.2).
	\begin{align*}
		\sum_{\substack{m,n \leq N \\ m \neq n}} \frac{a_m \bar{b}_n}{m - n} &\ll \Big(\sum_{n\leq N} |a_n|^2\Big)^{\half} 
		\Big(\sum_{n\leq N} |b_n|^2\Big)^{\half} \\
		\sum_{\substack{m,n \leq N \\ m \neq n}} \frac{a_m \bar{b}_n}{\log m - \log n} &\ll \sum_{n\leq N} n \left(|a_n|^2+|b_n|^2\right)
	\end{align*}
\end{prop}
\begin{prop}\label{simple:afe} 
For $a=a_R + i a_I \in \mathbb{C}$ and $2\pi S > |a_I|$, the zeta function may be approximated as 
(see Theorem 4.11 in Ref.\cite{titchmarsh1986theory})
	\begin{align}
	\zeta(a) = \sum_{n\leq S} n^{-a} + \frac{S^{1-a}}{a-1} + O(S^{-a_R})
	\end{align}
\end{prop}
\begin{prop}\label{chi_integral} 
Let $a,b$ be as in Theorem \ref{thm1}. Then for $y>0$ and $\gamma=\frac{1}{2}$
\begin{align}\label{chi-int}
	\int_{\gamma-iT}^{\gamma+iT} ds\, y^s\,\chi(s+a) \chi(1-s-b) &=  \frac{4\pi i }{1-a+b} \left(\frac{T}{2\pi}\right)^{1-a+b} \delta_{y,1}
	+ O\left(y^{1+\epsilon}T^{-\half+2 \epsilon}\right) \notag\\
	&\qquad + y^{1+\epsilon} |1\pm y|^{-1+2\epsilon} \Big(T^{2\epsilon} + |a_I-b_I|^{\frac{1}{2}+2\epsilon}\Big)
\end{align}
where $\delta_{y,1} = \begin{cases}
                                  1 &\text{if } y=1\\
                                  0 &\text{if } y\neq 1
                                  \end{cases}$, is the Kronecker delta.
\end{prop}
\noindent We discuss the proof of the above in some detail.
\begin{proof}
Recall the Mellin transformation of the cosine function
\begin{eqnarray}
	\int_{0}^{\infty} \frac{dx}{x} x^{s+b} \cos(2\pi x) = \frac{1}{(2\pi)^{s+b}} \cos\left(\tfrac{\pi}{2}(s+b)\right)\Gamma(s+b) = \frac{1}{2} \chi(1-s-b),
	\quad \Re(s+b) > 0
\end{eqnarray}
which allows us to replace one of the $\chi$-functions on the LHS of Eq.\eqref{chi-int}
\begin{align}\label{two-parts}
	\int_{\gamma-iT}^{\gamma+iT} ds\, y^s\,\chi(s+a) \chi(1-s-b)
	&= 2\int_{\gamma-iT}^{\gamma+iT}ds\, y^s \chi(s+a)\int_{0}^{\infty}\frac{dx}{x}x^{s+b}\cos(2\pi x)\notag\\
        = 2\int_{\gamma-iT}^{\gamma+iT}ds\, y^s \chi(s+a) & \left[\int_{0}^{\frac{T}{2\pi y}} \frac{dx}{x}x^{s+b}\cos(2\pi x) +
        \int_{\frac{T}{2\pi y}}^{\infty} \frac{dx}{x} x^{s+b} \cos(2\pi x)\right]
\end{align}
where we have split the inner integral into two parts. Let us evaluate the first term, which we begin by exchanging the order of integrations.
\begin{equation}\label{first-term}
	2 y^{1-a}\int_{0}^{\frac{T}{2\pi y}}dx\,x^{-a+b}\cos(2\pi x)\int_{\gamma-iT}^{\gamma+iT}ds(x y)^{s+a-1}\chi(s+a).
\end{equation}
and follow steps that are similar to those in the proof of Prop.~\ref{titch_lemma}, however, we shall need to interchange the integrals one more 
time.
	\begin{equation}\label{contour-terms}
	\int_{\gamma-iT}^{\gamma+iT} ds (x y)^{s+a-1} \chi(s+a) = \left( \oint_{\mathcal{C}}\; - \int_{\gamma+iT}^{b_{\gamma}+iT}\, 
	- \int_{b_{\gamma}+iT}^{b_{\gamma}-iT}\, - \int^{\gamma-iT}_{b_{\gamma}-iT} \right) ds\, (x y)^{s+a-1}\chi(s+a)
	\end{equation}
where $\mathcal{C}$ is a closed rectangular contour joining the points $\gamma-iT\to\gamma+iT \to b_{\gamma}+iT \to b_{\gamma}-iT \to 
\gamma-iT$ by horizontal or vertical straight lines. We take $b_{\gamma}=T^P$ for a suitably large power $P>0$ and assume that $|a_I |< T$. 
The integral over the closed contour $\mathcal{C}$ gets contributions from the residues at the poles of the integrand, which are at the poles of
the gamma function in $\chi$.
\begin{align}\label{residue-term}
	2\oint_{\mathcal{C}}ds\, (2\pi x y)^{s+a-1} \Gamma(1-s-a) \sin\left(\tfrac{\pi}{2}(s+a)\right) &=
	4\pi i\, \sum_{l=0}^{\lfloor {b_\gamma}/{2}\rfloor} \frac{(-1)^{l}}{(2l){!}} (2\pi x y)^{2l} \notag\\
	&= 4 \pi i\, \cos(2 \pi x y) - R, 
	\end{align}
where $R$ is the remainder of the Taylor series expansion of the cosine function truncated at the $2 \lfloor b_\gamma/2\rfloor$-th term. 
The bound $|R| \le {|2 \pi x y|^{2\lfloor \frac{b_\gamma}{2}\rfloor+1}}/{\left(2\lfloor b_\gamma/2\rfloor + 1\right)!}$ for $xy \le T/(2\pi)$ and 
$b_{\gamma} = T^P$ for $P > 1$, implies that $R$ is exponentially small. Hence its contribution to Eq.\eqref{first-term} is $o(1)$. 
	
Let us now estimate the remaining parts of the rectangular contour in Eq.\eqref{contour-terms}. The contribution from the integrals on the two 
horizontal contours are the same, therefore, we shall discuss the details of the first one. We change variable in
\begin{equation}\label{gamma-to-bgamma}
	\mathcal{I}_{H_1} = -2 \int_{\gamma+iT}^{b_{\gamma}+iT} ds\, (2\pi x y)^{s+a-1} \sin\left(\tfrac{\pi}{2}(s+a)\right)\, \Gamma(1-s-a)
\end{equation}
to integrate in $\sigma = \Re s$ from $\gamma$ to $b_\gamma$. Next we use Stirling's approximation \cite{andrews1999special}
\begin{equation}\label{eq:Stirling}
\Gamma(z) \stackrel{|z|\to\infty}{\longrightarrow} \sqrt{2\pi}\, z^{z-\frac{1}{2}} e^{-z},\:\text{ for } |\text{arg}(z)| \le \pi - \delta, \;\delta\to 0+
\end{equation}
to find a bound on the integrand
\begin{equation}\label{horiz-term}
\left| \mathcal{I}_{H_1} \right| \lesssim \frac{\sqrt{2\pi}}{2\pi xy} \int_\gamma^{b_\gamma} d\sigma\, \sqrt{H(\sigma)}\, 
\exp \left[ - (\sigma + a_R) \log \frac{H(\sigma)}{2\pi e xy} - (T + a_I)\left(\frac{\pi}{2} - \theta(\sigma)\right)\right] 
\end{equation}
where $H(\sigma) = \sqrt{(\sigma + a_R)^2 + (T + a_I)^2}$ and $\tan\theta(\sigma) = (T+a_I)/(\sigma+a_R)$. The integral over the other 
horizontal segment in Eq.\eqref{contour-terms} from $b_{\gamma}-iT\to\gamma-iT$ results in a similar estimate.
We note that the argument of the exponential in Eq.\eqref{horiz-term}, which appears in the integral over $\sigma$, is a monotonically 
decreasing function whose derivative is $ - \log \frac{H(\sigma)}{2\pi x y}$, which decays fast. Therefore, $\mathcal{I}_{H_1} \ll \frac{1}{x y}
(x y)^{\gamma+a_R} T^{a_R}$ and 
\begin{equation}\label{horiz-estimate}
	O\left(y^{\half} \left(\frac{T}{y}\right)^{b_R - \frac{1}{2}}\right) = O\left(y^{1+\epsilon} T^{-\frac{1}{2}+\epsilon}\right)
\end{equation}
in Eq.\eqref{first-term}. Here we used the estimate $\int_0^A x^{\alpha} \cos(2\pi x) = O(A^{\alpha})$ by integration by parts. 

The estimation of the vertical integral from $b_{\gamma}-iT\to b_{\gamma}+iT$ in Eq.\eqref{contour-terms} is done in exactly the same 
fashion. This time we change variable to 
$\tau = \Im{s}$ in $\mathcal{I}_V = - \displaystyle{\int^{b_\gamma-iT}_{b_{\gamma}+iT}} ds\, (x y)^{s+a-1}\chi(s+a)$ 
to get an integral from $-T$ to $T$. Once again we use Stirling's approximation to find 
\begin{equation*}
\left| \mathcal{I}_V\right| \ll \frac{1}{xy}\sqrt{\frac{2}{\pi}} \int_0^T d\tau\, \sqrt{V(\tau)} \exp\left[ - (b_\gamma + a_R)\log \left(
\frac{V(\tau)}{2\pi exy}\right) - (\tau + a_I)\left(\frac{\pi}{2}-\varphi(\tau)\right)\right]
\end{equation*}
where $V(\tau) = \sqrt{(b_\gamma+a_R)^2 + (\tau+a_I)^2}$ and $\tan\varphi(\tau) = (\tau+a_I)/(b_\gamma+a_R)$. 
Since $b_{\gamma} = T^{P}$ with $P>1$, this falls off exponentially in $T$ and thus its contribution to Eq.\eqref{first-term} is 
exponentially smaller than Eq.\eqref{horiz-estimate}.
	
We shall now evaluate the contribution from the dominant term in Eq.\eqref{contour-terms} to Eq.\eqref{first-term} after performing the 
integral over $x$.
\begin{align*}
	 \frac{8\pi i}{y^{a-1}} \int_{0}^{\frac{T}{2\pi y}} dx\, x^{b-a} \cos(2\pi x) \cos(2\pi x y)
	 = \frac{4\pi i}{y^{a-1}} \int_{0}^{\frac{T}{2\pi y}} dx\, x^{b-a} \left(\cos(2\pi x(1-y)) + \cos(2\pi x(1+y))\right)
\end{align*}
When $y=1$ in the first term above, a simple integration yields 
\begin{align*}
	\frac{4 \pi i}{b-a+1}  \left(\frac{T}{2\pi }\right)^{b-a+1}
\end{align*} 
which is the main term of our interest. For the remaining integrals, we note that 
\begin{align}\label{cos-integral}
	\frac{4\pi i}{y^{a-1}} \int_{0}^{\frac{T}{2\pi y}} & dx\, x^{b-a} \cos\left(2\pi x(1\pm y)\right)\: =\: 
	\frac{4 \pi i}{y^{b}} \Big( \frac{T}{2\pi }\Big)^{b-a+1} \int_{0}^{1} d\omega \, \omega^{b-a} \cos\left(\frac{(1\pm y)}{y}T\omega\right)\notag\\
	&= \frac{2\pi i}{y^b(b-a+1)}\bigg(\frac{T}{2\pi}\bigg)^{b-a+1} \left[\,{}_{1}F_{1}\left(b-a+1, b-a+2 ; i\frac{(1\pm y)}{y}T\right)\right.\notag\\
	& \qquad\qquad\qquad\left. +\, {}_{1}F_{1}\left(b-a+1, b-a+2; -i \frac{(1\pm y)}{y}T\right)\right]
\end{align}
We now use the asymptotic form of the confluent hypergeometric function 
\begin{equation}\label{confl-hyp-fn}
	{}_{1}F_{1} (a,b; z) \stackrel{|z| \to \infty}{\longrightarrow} \Gamma(b) \left( \frac{e^{z}z^{a-b}}{\Gamma(a)} + 
	\frac{(-z)^{-a}}{\Gamma(b-a)}\right)
\end{equation}
to put a bound on all terms (except the first term with $y=1$) in Eq.\eqref{cos-integral}. The contributions of these are
\begin{align}\label{cosint1}
	\text{LHS}&(\text{Eq.}\eqref{cos-integral})\Big|_{(1-y)\ne 0}\\
	& \ll \frac{T^{b-a+1}}{y^b} \left[ \frac{\sin\left(\frac{(1\pm y)}{y}T\right)}{\frac{(1\pm y)}{y}T} + \left(\frac{(1\pm y)}{y}T\right)^{a-b-1}
	\Gamma(b-a+1) \cos\left(\tfrac{\pi}{2}(b-a+1)\right) \right] \notag\\
	&\ll y^{1+\epsilon} \left|1\pm y\right|^{-1+2\epsilon} \left( T^{2\epsilon }+ |a_I-b_I|^{\frac{1}{2}+2\epsilon} \right)
\end{align}
In the above, we have also used the result
\begin{equation}\label{chi_asymptotic}
\chi(\sigma + it) 
=\left| \frac{t}{2\pi} \right|^{\frac{1}{2} - \sigma - it}e^{it+\frac{i\pi}{4} \text{sign}(t)} \left(1 + O \left(\frac{1}{|t|}\right)\right)
\end{equation}
from Ref.\cite{titchmarsh1986theory} (see Eq.(4.12.3)).

Finally, it remains to estimate the second term of \eqref{two-parts}. To this end, we first rescale the variable of integration so that we may
use Eq.\eqref{cos-integral} and the asymptotic form of the confluent hypergeometric functions Eq.\eqref{confl-hyp-fn} in the limit of large 
values of $T$, to write 
\begin{eqnarray}\label{final-int}
	&&\!\!\!\!\!\!\!\!\!\!\!\!\!\!\!\!\!\!\!\!\!\!\!\!\!\!\!\!\!\!\!\!
	2 \int_{\gamma-iT}^{\gamma+iT} ds\, y^s \chi(s+a) \int_{\frac{T}{2\pi y}}^{\infty} \frac{dx}{x}x^{s+b} \cos(2\pi x)\notag\\
	&=& \frac{2}{y^b} \int_{\gamma-iT}^{\gamma+iT}ds\, \chi(s+a) \bigg(\frac{T}{2\pi}\bigg)^{s+b}
	\left[ \frac{\pi (y/T)^{s+b}}{2\Gamma(1-s-b) \sin\left(\frac{\pi}{2}(s+b)\right)} \right.\notag\\
	&& \qquad \left. -\, \frac{1}{2(s+b)} \left({}_{1}F_{1}(s+b, s+b+1; \frac{iT}{y}\right) + {}_{1}F_{1}\left(s+b, s+b+1 ; -\frac{iT}{y}\right)
	\right]  \notag\\
	&\stackrel{T\to\infty}{\longrightarrow}& 
	\frac{2}{y^b} \int_{\gamma-iT}^{\gamma+iT} ds\, \chi(s+a) \bigg(\frac{T}{2\pi}\bigg)^{s+b} 
	\left[ \Gamma(s+b)\cos\left(\tfrac{\pi}{2}(s+b)\right) \left(\frac{y}{T}\right)^{s+b} \right. \notag\\
	&&\qquad\qquad -\, \left. \frac{\sin(T/y)}{T/y} - \Gamma(s+b) \cos\left(\tfrac{\pi}{2}(s+b)\right) \left(\frac{y}{T}\right)^{s+b} \right] \notag\\
	&\ll& \sin(\tfrac{T}{y}) \frac{T^{b-a}}{y^{b-1}} \int_{\gamma-iT}^{\gamma+iT} ds\, \chi(s+a) \left(\frac{T}{2\pi}\right)^{s+a-1} 
	\:\sim\: O\left( y^{1+\epsilon}T^{-\half+2 \epsilon}\right)
\end{eqnarray}
In the last step we use Lemma \ref{titch_lemma} with $x = \frac{T}{2 \pi}$ and $\gamma = \frac{1}{2} + a_R$.
\end{proof} 

\section{Computation of the second moment}\label{sec:CompMom2}
We now compute the second moment. Let us rewrite $I(a,b,T)$ as a contour integral
\begin{align}\label{afe-product}
    I(a,b,T) &= \int_{-T}^T dt\, \zeta(\tfrac{1}{2}+a+it)\zeta(\tfrac{1}{2}-b-it)\notag\\
		&= -i \int_{\half-iT}^{\half+iT} ds\, \zeta(s+a)\zeta(1-s-b)\notag\\
	 	&= -i \int_{\half-iT}^{\half+iT} ds\, \Big(\sum_{n\le q} n^{-s-a} + \chi(s+a) \sum_{n\le q} n^{s+a-1} + E(s+a)\Big)\notag\\
		&\qquad \times\,\Big(\sum_{m\le q'}m^{s+b-1}+\chi(1-s-b)\sum_{m\le q'}m^{-s-b}+E(1-s-b)\Big)
\end{align}
in which we have substituted the approximate functional equation Eq.\eqref{afe} in the integrand with 
$q = \left\lfloor\sqrt{\tfrac{|t+a_I|}{2\pi}}\right\rfloor$ for $\zeta(s+a)$ and 
$q' = \left\lfloor\sqrt{\frac{|t+b_I|}{2\pi}}\right\rfloor$ for $\zeta(1-s-b)$. 

First, we look at the diagonal terms of Eq.\eqref{afe-product} with $z = \left\lfloor\sqrt{\tfrac{1}{2\pi}\min(|t+a_I|,|t+b_I|)}\right\rfloor$.
\begin{align*}
-i \int_{\half-iT}^{\half+iT} ds\, & \Big(\sum_{n\le z} n^{b-a-1} + \chi(s+a)\chi(1-s-b) \sum_{n\le z} n^{a-b-1}\Big)\\
&= -i \int_{\half-iT}^{\half+iT} ds\, \Big[ \zeta(1+a-b) + \frac{z^{b-a}}{b-a} + O(z^{b_R-a_R-1})\\
&\qquad\quad+\; \psi(a,b,t) \Big(\zeta(1-a+b) + \frac{z^{a-b}}{a-b} + O(z^{a_R-b_R-1})\Big) \Big] 
\end{align*}
where we adapted Prop.~\ref{simple:afe} as appropriate. This correctly reproduces the main terms of Eq.\eqref{ingham}, while the 
$O$-terms give the order $O(T^{\half+\epsilon})$. As for the remaining terms, we use  
\begin{align*}
	\int_{\half-iT}^{\half+iT}ds\, &\psi(a,b,t)\, \frac{z^{a-b}}{a-b} = \frac{4\pi i}{(1-a+b)(a-b)}\Big[\left(\frac{T}{2\pi}\right)^{1+\frac{b-a}{2}} + 
	O(T^{3\epsilon}) +O(|a_I-b_I|^{\half+2\epsilon}T^{\epsilon})\\
	&- \int_{\half-iT}^{\half+iT}ds \left(\frac{t}{2\pi}\right)^{1+b-a} \frac{a-b}{2}\left(\frac{t}{2\pi}\right)^{-1+\frac{a-b}{2}}\Big]\\
	&=\frac{4\pi i}{(1-a+b)(a-b)} \left(\frac{T}{2\pi}\right)^{1+\frac{b-a}{2}} \Big[ 1-\frac{a-b}{2(1+ \frac{b-a}{2})}\Big] + O(T^{3\epsilon})+
	O(|a_I-b_I|^{\half+2\epsilon}T^{\epsilon})\\
	 &=  \int_{\half-iT}^{\half+iT} ds\, \frac{z^{b-a}}{a-b} + O(T^{3\epsilon})+O(|a_I-b_I|^{\half+2\epsilon}T^{\epsilon})
\end{align*}
which is obtained from Prop. \ref{chi_integral} and integration by parts.

Now we come to the computation of the non-diagonal terms of Eq.\eqref{afe-product}.
\begin{itemize}
\item
 First, we consider the terms which do not involve the $\chi$-function.
\begin{align}\label{non-diag1}
	-i \sum_{\stackrel{m,n\ll \sqrt{T}}{n\ne m}} \frac{m^{b-1}}{n^a} \int_{\half-iT}^{\half+iT} ds\, \frac{m^{s}}{n^{s}} 
	&= \sum_{\stackrel{m,n\ll \sqrt{T}}{n\ne m}} \frac{m^{b-\half}}{n^{a+\half}}\, \frac{(m/n)^{iT} - (m/n)^{-iT}}{\log(m/n)}\\
	&= O \Big(\sum_{n\leq N} n \left( n^{-2 a_R - 1} + n^{ 2 b_R - 1}\right)\Big) \\
	&= O(T^{\half-a_R})+O( T^{\half+b_R}) \:=\: O(T^{\half+\epsilon})
\end{align}
In the above, we have estimated the sum in the first line by Prop. \ref{hilbert} with $a_n = n^{-\half+b\pm iT}$, 
$b_n = n^{-\half-\overline{a}\pm iT}$ and $N=O(\sqrt{T})$.
\item
Next we consider the terms involving products of two $\chi$'s.
\begin{align}\label{non-diag2-estimate1}
	-i \sum_{\stackrel{m,n\ll \sqrt{T}}{n\ne m}} \frac{n^{a-1}}{m^b} 
	& \int_{\half-iT}^{\half+iT} ds\, \frac{n^{s}}{m^{s}}\, \chi(s+a)\chi(1-s-b)\\
	&\ll \sum_{\stackrel{m,n\ll \sqrt{T}}{n\ne m}} \frac{n^{a-1}}{m^b}\, \left[\Big( \frac{n}{m} \Big)^{1+\epsilon}
	T^{-\half +2\epsilon} + \frac{n^{1+\epsilon} m^{-3\epsilon}}{(m-n)^{1 - 2\epsilon}}
	\left(T^{2\epsilon} + |a_I - b_I|^{\tfrac{1}{2}+2\epsilon}\right)\right]\\
	&\ll T^{\half+\epsilon} + \sum_{\stackrel{m,n\ll \sqrt{T}}{n\ne m}} \frac{n^{a+\epsilon}m^{-b-\epsilon}}{m-n} \left(T^{2\epsilon} + 
	|a_I-b_I|^{\frac{1}{2}+2\epsilon}\right)\\
	&= O\Big( T^{\half(1+ a_R-b_R)} \big(T^{2\epsilon} + |a_I-b_I|^{\frac{1}{2} + 2\epsilon}\big)\Big)\\ 
	&= O(T^{\half+3\epsilon}) + O\big(T^{\half+\epsilon}|a_I-b_I|^{\frac{1}{2}+2\epsilon}\big)
\end{align}
In this case, we have used Prop. \ref{chi_integral} with $y=n/m$ at the first step, followed by an estimation of the second term by using Prop. 
\ref{hilbert}.
\item
Now let us evaluate the two cross terms with one factor of the $\chi$-function. One of these is estimated as 
\begin{align}\label{cross-term1}
-i\int_{\half-iT}^{\half+iT} ds\, & \sum_{n \le \sqrt{\tfrac{|t+a_I|}{2\pi}}}\; \sum_{m \le \sqrt{\tfrac{|t+b_I|}{2\pi}}} 
\frac{1}{n^{s+a} m^{s+b}}\, \chi(1-s-b)\\
&\ll \sum_{m,n\ll \sqrt{T}} n^{b-a} \displaystyle{\int_{\half+i \max(2\pi n^2+a_I,2\pi m^2+b_I)}^{\half+iT}} 
\frac{ds}{(nm)^{s+b}}\, \chi(1-s-b)\\
&= \sum_{m,n\ll \sqrt{T}} n^{b-a} \Big[ O\Big( \frac{T^{b_R}}{(nm)^{\half+b_R} \log \tfrac{T}{2\pi n m}} \Big) + 
O\Big( \frac{T^{b_R}}{(nm)^{\half+b_R}} \log T\Big)\Big] 
\end{align}
using Prop. \ref{titch_lemma} with $x= nm$, $\gamma = \half+b_R$ and both limits of the integral ($T$ as well as $2\pi \max(n^2,m^2)$) 
for $T$ to cancel the leading term $2 \cos(2\pi m n)$. As for the other term 
\begin{align} \label{cross-term2}
-i\int_{\half-iT}^{\half+iT} ds & \sum_{n \le \sqrt{\tfrac{|t+a_I|}{2\pi}}}\; \sum_{m \le \sqrt{\tfrac{|t+b_I|}{2\pi}}} n^{s+a-1} m^{s+b-1}\, \chi(s+a)\\
&= \sum_{m,n\ll \sqrt{T}} m^{b-a} \Big[ O \Big( \frac{(nm)^{a_R - \half}}{T^{a_R} \log\tfrac{T}{2\pi n m}} \Big) +
O\Big( \frac{(nm)^{a_R - \half}}{T^{a_R}} \log T \Big) \bigg] 
\end{align}
where, similar to previous case, we have changed the variable of integration to $1-s-a$ and uses Prop. \ref{titch_lemma} with $x= nm$, 
$\gamma = \half-a_R$ and both the limits $T$ as well as $2\pi \max(n^2,m^2)$ for $T$.

In order to complete the estimates, we analyse the  sums in Eqs.\eqref{cross-term1} and \eqref{cross-term2}. The first error terms are 
subdivided into the ranges $nm \leq \tfrac{T}{4\pi}$ and $\tfrac{T}{4\pi} < nm < \tfrac{T}{2\pi}$, while the second  error terms are handled 
like the sum in the first range. 

In Eq.\eqref{cross-term1} the sum in the range $nm\leq T/4 \pi$ is estimated to be
\begin{equation}\label{cross-term1-estimate1}
	O\Big(\, T^{b_R}\sum_{m,n\ll \sqrt{T}} n^{-\half-a_R} m^{-\half-b_R} \Big)
	= O\left(T^{\half+\frac{b_R-a_R}{2}}\right) = O(T^{\half+\epsilon})
\end{equation}
The sum in \eqref{cross-term1} in the range 
$\tfrac{T}{4\pi} < nm < \tfrac{T}{2\pi}$ can be estimated as
\begin{align}\label{cross-term1-estimate2}
	O\bigg( \sum_{\frac{T}{4\pi} < k < \frac{T}{2 \pi}}\, \sum_{\substack{m\ll \sqrt{T}\\ m | k}}   
	\frac{m^{a_R-b_R}\, T^{1+b_R}}{T^{\half+b_R} \left(\frac{T}{2\pi} - k\right)} \bigg)
	&= O\Big( T^{\frac{1}{2} (a_R-b_R+1)} \sum_{\frac{T}{4\pi} < k < \frac{T}{2 \pi}} \frac{d(k)}{\frac{T}{2\pi} - k} \Big)\notag\\ 
	& \!\!\!\!\!\!\!\!\!\! = O\big(T^{\half(a_R-b_R+1)+\epsilon}\log T\big) \:=\: O(T^{\half+2\epsilon})
\end{align}
where $d(k)$ is the divisor function for which we have used $d(k) = O(k^{\epsilon})$(see \cite{titchmarsh1986theory}, p.145 or 
\cite{ivic2013riemann}, Eq. (1.7.1)). The sum in Eq.\eqref{cross-term2} in both the ranges, can be estimated analogously.
\item 
Finally we estimate the contributions from the products involving the error terms $E(s+a)$ and $E(1-s-b)$ in Eq.\eqref{afe-product}. The 
convexity bound $\zeta(\sigma+it) \ll |t|^{\half(1 - \sigma)+\delta}$ for $0< \sigma <1$ and any $\delta>0$ (see Ref.\cite{titchmarsh1986theory}, 
Eq.(5.1.4)), implies that each of the Dirichlet polynomials in Eq.\eqref{afe-product} can be estimated to be $O\left( |\zeta(s+a)| + |\zeta(1-s-b)| 
\right)= O\left(|t|^{\frac{1}{4}+\epsilon}\right)$. Therefore, the contribution of $O(t^{-3/4})$ from \eqref{siegel-term} multiplied with any of the 
Dirichlet polynomials is
\begin{align}\label{error-estimate}
	O(T^{1+\frac{1}{4}-\frac{3}{4}+\epsilon}) = O(T^{\half+\epsilon}).
\end{align} 
Now we consider contribution of the terms involving $\phi(t+a_I)$ and $\phi(t+b_I)$ from Eq.\eqref{siegel-term}, where
\begin{align*}
	\phi(t) = \frac{\cos\left(t-(2q+1)\sqrt{2\pi t}-\frac{\pi}{8}\right)}{\cos \sqrt{2\pi t}}	
	= \frac{\cos\left(2 \pi (f^2-f)-\frac{\pi}{8}\right)}{\cos 2\pi f}
\end{align*} 
where $q=\left\lfloor \sqrt{t/2 \pi}\right\rfloor$ and $f=\sqrt{t/2 \pi}-q $ is the fractional part of $\sqrt{t/2 \pi}$. Thus, $\phi$ is a periodic and 
continuous function of $\sqrt{t/2 \pi}$. Moreover, $(-1)^{q-1}\phi(t+a_I)$ is bounded and monotonic in the intervals  $\left(2 \pi k^2, 2 \pi 
(k+\half)^2\right)$ and $\left(2 \pi (k+\half)^2, 2 \pi (k+1)^2\right)$ for integers $k$. 
The product of the term involving $\phi(t+a_i)$, with two different Dirichlet polynomials from \eqref{afe-product} needs to be estimated. 
One of these products is 
\begin{align*} 
	&\int_{-T}^{T}\! dt \!\!\!\!\!\sum_{m \leq \sqrt{\frac{|t+b_I|}{2\pi}}} m^{-\half+b_R+i(t+b_I)} e^{- \frac{i\pi}{2}(s-1) - \frac{i}{2} (t+a_I) - 
	\frac{i\pi}{8}} (2\pi |t+a_I|)^{\frac{1}{2}(s-1)} \Gamma(1-s) (-1)^{q-1}\phi(|t+a_I|) \nonumber \\
	&=\sum_{m \ll \sqrt{T}} \int_{I_m} dt\, e^{- \frac{i\pi}{2}(s-1) - \frac{i}{2} t - \frac{i\pi}{8}} (2\pi |t+a_I|)^{\frac{1}{2}(s-1)} 
	\Gamma(1-s) m^{-\half+b_R+i(t+b_I)}(-1)^{q-1}\phi(|t+a_I|) 
\end{align*}
where $s=\half+i(t+a_I)$ and $I_m = [-T, - 2\pi m^2-b_I] \cup [2\pi m^2-b_I,T]$. Using Stirling's approximation \eqref{eq:Stirling} and the 
substitution $t+a_I\rightarrow t$, we get
\begin{align}\label{siegel-integral2}	
	&\ll \sum_{m \ll \sqrt{T}}m^{-\half+b_R} \int_{J_m} dt\, |t|^{-\frac{1}{4}-\frac{a_R}{2}}\,  (-1)^{q-1}\phi(|t|) \exp(i F_m(t))\
\end{align}
where $F_m(t) = - \frac{t}{2}\log(\frac{t}{2\pi})+\frac{t}{2} + (t+b_I-a_I)\log m$ and the domain of integral is 
$J_m=I_m+a_I = [-T+a_I, - 2\pi m^2+a_I-b_I] \cup [2\pi m^2+a_I-b_I,T+a_I]$. We note that $F_m'(t) = -\frac{1}{2}\log(\frac{t}{2\pi})+ \log m$ is 
monotonic and nonvanishing for $t \neq 2\pi m^2$. Hence the function $|t|^{-\frac{1}{4}-\frac{a_R}{2}}(-1)^{q-1} \phi(|t|)/F_m'(t)$ is monotonic 
and bounded over intervals $(2\pi k^2,2\pi (k+\half)^2)$ and $(2\pi (k+\half)^2,2\pi (k+1)^2)$ for integers $k\neq m,m-1$, where 
Lemma 4.3 of Titchmarsh \cite{titchmarsh1986theory} can be used to estimate the integral in \eqref{siegel-integral2} as 
\begin{equation*}
	\frac{(k^2)^{-\frac{1}{4}-\frac{a_R}{2}}}{\log k - \log m}
\end{equation*}	 
Thus the overall estimate for Eq.\eqref{siegel-integral2} over the region $J_m \setminus [2\pi (m-1)^2,2\pi (m+1)^2]$ is 
\begin{equation}\label{siegel-sum}
 \sum_{m \ll \sqrt{T}} \sum_{\substack{k \ll \sqrt{T} \\ k \neq m-1,m}} \frac{m^{-\frac{1}{2}+b_R}k^{-\frac{1}{2}-a_R}}{\log k - \log m} \ll 
 \sum_{k \ll \sqrt{T}} k(k^{-1+2 b_R}+k^{-1-2 a_R} )	= O\left(T^{\frac{1}{2}+\epsilon}\right)
\end{equation}
where we have used Proposition \ref{hilbert} with $N=O(\sqrt{T})$, $a_m = m^{-\frac{1}{2}+b_R}$ and $b_n = n^{-\frac{1}{2}-a_R}$.

For the region $J_m \cap [2\pi (m-1)^2,2\pi (m+1)^2]$, we use the trivial estimate for the integral in Eq.\eqref{siegel-integral2} as $m^{-\half-a_R} m$. Hence the overall estimate for Eq. \eqref{siegel-integral2} in this region is
\begin{equation}\label{siegel-sum2}
	\sum_{m \ll \sqrt{T}} m^{-\frac{1}{2}+b_R-\half-a_R+1}= O\left(T^{\frac{1}{2}+\epsilon}\right)
\end{equation}

The product of the error term involving $\phi(t+a_i)$, with the other Dirichlet polynomial from Eq.\eqref{siegel-term} is
\begin{align}\label{siegel-integral3}
	&\int_{-T}^{T} dt\: E(s+a) \chi(1-s-b) \sum_{m \leq \sqrt{\frac{|t+b_I|}{2\pi}}} m^{-\half-b_R-i(t+b_I)} \nonumber \\
	&=\sum_{m \ll \sqrt{T}} \int_{I_m} dt\: E(s+a) \chi(1-s-b) m^{-\half-b_R-i(t+b_I)}\\
	&\ll \sum_{m \ll \sqrt{T}} m^{-\half-b_R} \int_{I_m} dt\, |t+b_I|^{b_R}|t+a_I|^{-\frac{1}{4}-\frac{a_R}{2}}(-1)^{q-1}\phi(|t+a_I|) \exp(i G_m(t+a_I))
\end{align}
where $s=\half+it$ and $G_m(t)= -\frac{t}{2}- \half t \log(\frac{t}{2\pi})+ (t+b_I-a_I)\log(\frac{t+b_I-a_I}{2\pi})-(t+b_I-a_I)\log m$. Here we have 
used Eq.\eqref{chi_asymptotic} to get an asymptotic for $\chi(1-s-b)$ and Stirling's approximation \eqref{eq:Stirling} for the rest.

Now we note that $G_m'(t)=-\half \log \frac{t}{2 \pi}+ \log \frac{t+c_I}{2 \pi}-\log m$ , where $c_I=b_I-a_I$. We first consider the case when 
$\pi m^2\geq 2 c_I$. The critical points $G_m'(t)=0$ occur at
\begin{equation*}
	t+c_I = \pi m^2 \pm \sqrt{\pi m^2(\pi m^2-2 c_I)}
\end{equation*}
leads to the critical point $m^2 \ll t+c_I \ll m^2$ and this case therefore can be handled completely analogous to the estimation 
for \eqref{siegel-integral2}.

For the other case $\pi m^2< 2 c_I$, we get that $G_m(t)$ is monotonic in the entire region $I_m$ with the minimum value for $G_m'(t)$ being 
$\half \log \frac{2 c_I}{\pi m^2}$. Hence the contribution to \eqref{siegel-integral3} from these terms is
\begin{equation}\label{siegel-sum3}
	\sum_{m \ll \sqrt{c_I}}\frac{ m^{-\frac{1}{2}+b_R}}{\log(2 c_I)-\log(\pi m^2)} \sum_{k\ll \sqrt{T}} k^{-\frac{1}{2}-a_R+2b_R}
	= O\left(T^{\frac{1}{4}+3\epsilon} |a_I-b_I|^{\frac{1}{4}+\epsilon}\right)
\end{equation}
\end{itemize}

\section{Shifted second moment of the zeta function as a kernel}\label{sec:ZKernel}
In this section we present the proof of Theorem \ref{thm2}.  
\begin{proof} 
We start with Theorem \ref{thm1} with $\text{max}(\alpha,\beta) = O(T)$.
\begin{align}\label{eq:Ingham-Bettin}
	\int_{-T}^{T}dt\, \zeta\left(\tfrac{1}{2}+i\alpha+it\right)\zeta\left(\tfrac{1}{2}-i\beta-it\right)&=\int_{-T}^{T}d\tau\,\left(\zeta(1+i\alpha-i\beta) + 
	\psi(i\alpha,i\beta,\tau)\zeta(1+i\beta-i\alpha)\right)\notag\\
	&\qquad +\, O(T^{\frac{1}{2}+3 \epsilon})+O(T^{\frac{1}{2}}|\alpha-\beta|^{\frac{1}{2}+2\epsilon})\notag\\
	&= 2T\,\zeta(1+i\alpha-i\beta) + 2T\, \Big[\frac{\left({T}/{2\pi}\right)^{i\beta-i\alpha}}{1-i\alpha+i\beta} \notag\\
	&\qquad +\, O(T^{-\frac{1}{2}+2\epsilon}) + 
	O(T^{-1}|\alpha-\beta|^{\frac{1}{2}+2\epsilon})\Big] \zeta(1+i\beta-i\alpha)\notag\\
	&\qquad +\, O(T^{\frac{1}{2}+3 \epsilon}) + O(T^{\frac{1}{2}} |\alpha-\beta|^{\frac{1}{2}+2\epsilon})
\end{align}
We have applied Proposition \ref{chi_integral} to arrive at the expression above. Now we use \eqref{eq:Ingham-Bettin} to evaluate the required 
integral for fixed $\beta$, with $\text{max}(\beta) < \text{max}(\alpha) \sim R \sim T$.
\begin{align}\label{eq:Kernellimt}
	\lim_{T\to\infty} \int_{-R}^{R} d\alpha\, f(\alpha)  
	&\left[ \frac{1}{2T} \int_{-T}^{T} \zeta\left(\tfrac{1}{2}+it+i\alpha\right) \zeta\left(\tfrac{1}{2}-it-i\beta\right)\, dt\right]\notag\\
	&=  \lim_{T\to\infty} \int_{-R}^{R} d\alpha\, f(\alpha) \Big(\zeta(1+i\alpha-i\beta) + \left(\frac{T}{2\pi}\right)^{i\beta-i\alpha}
	\frac{\zeta(1+i\beta-i\alpha)}{1+i\beta-i\alpha}\notag\\ 
	&\qquad + O(T^{-\frac{1}{2}+3\epsilon})+O(T^{-\frac{1}{2}} |\alpha-\beta|^{\frac{1}{2}+2\epsilon})\Big)
\end{align}
The test function $f(\alpha)$ satisfies the conditions: $\displaystyle{\int_{-\infty}^{\infty} d\alpha\, |\alpha|^{\half+2\epsilon}\, \left|f(\alpha)\right|}$ 
and $\displaystyle{\int_{-\infty}^{\infty} d\alpha\, \left|f(\alpha)\right|}$ exist, therefore, the last two terms in the integral in \eqref{eq:Kernellimt} 
tend to $0$ as $T\rightarrow\infty$. We evaluate the remaining two terms in the integral in \eqref{eq:Kernellimt} along a closed contour 
$\mathcal{C}$ in the complex $\alpha$-plane consisting of the straight segment along the real axis from $-R$ to $\beta - \epsilon$, followed by 
a semi-circle around $\beta$ (where the $\zeta$-function has a simple pole) of infinitesimally small radius $\epsilon$, continuing along the real 
axis from $\beta +\epsilon$ to $R$, which is then closed along a semi-circular Jordan curve $\Gamma_R$ connecting $R$ to $-R$ on the upper 
or the lower half plane, as appropriate. Note that the integrand is bounded everywhere. Further, in the limit $\alpha\to\beta$, the singularity 
coming from the two terms cancel. Hence, it is sufficient to compute the principal value (PV) of each term and add them. 

Let us consider the contribution from the Jordan curve to the first term. On $\Gamma_R$, we parametrise $\alpha = \beta + Re^{i\theta}$, 
and use the reflection identity of the zeta function to write 
$\zeta(1 + i \alpha-i\beta) = 2 (2\pi)^{iRe^{i\theta}} \mathrm{cosh}\left(\frac{\pi}{2}R e^{i\theta}\right) \Gamma(-iRe^{i\theta}) \zeta(-iRe^{i\theta})$. 
The requirement that $|\zeta(R\sin\theta-iR\cos\theta)|$ is bounded dictates that we close the contour on the UHP, which is consistent with the 
condition on the argument of $\Gamma(-iRe^{i\theta})$ to determine its asymptotic behaviour. However, the contribution from the 
$\Gamma$-function makes $\left|\zeta(1+i\alpha-i\beta)\right|$ grow as $R^{R\sin\theta}$. Since the test function $f(\alpha)= o(R^{-R})$, the 
integral along the Jordan curve is $o(1)$, as $R\rightarrow\infty$. Thus, we find
\begin{equation*}
	 \mathrm{PV}\int_{-R}^{R} d\alpha\, f(\alpha) \zeta(1+i\alpha-i\beta) -\pi f(\beta) + o(1)
	= 2\pi i\!\!\!\sum_{\alpha_A\in\text{UHP}} \!\!\mathrm{Res}_{\alpha_A}\left(f(\alpha_A) \zeta(1+i\alpha_A-i\beta)\right) = 0
\end{equation*}
where the sum is over the poles $\alpha_A$ of $f$ within the contour in the $\text{UHP}$, above the pole $\alpha=\beta$ from the zeta function. 
This is zero as $f(\alpha)$ has been chosen to be analytic there. Hence, we obtain
\begin{align}\label{part1}
	\mathrm{PV}\int_{-R}^{R} d\alpha\, f(\alpha) \zeta(1+i\alpha-i\beta)  = \pi f(\beta) + o(1)
\end{align}
We follow the same steps for the second integral. The semi-circular contour is to be closed on the LHP to ensure that the zeta function (after
using the reflection identity) remains bounded. This is also consistent with the condition on the asymptotics of the $\Gamma$-function. Taking 
$\alpha = \beta + Re^{i\theta}$ and $R \sim T$, the leading order asymptotics from the integrand is $|f(\beta + Re^{i\theta}) R^{-\frac{3}{2}}
e^{g(\theta)R}|$, where $g(\theta)$ is a bounded function of $\theta$. Thus it would be sufficient to take $f(Re^{i\theta}) \sim e^{-cR}$, for 
some positive $c > g(\theta)$ for $\theta \in (-\pi, 0)$. The integral from the Jordan curve is $O(R^{-\delta})$ for some $\delta>0$.

Therefore, the contribution from the second integral is
\begin{align*}
	\mathrm{PV}\int_{-R}^{R} d\alpha\, & f(\alpha) \left(\frac{T}{2\pi}\right)^{i\beta-i\alpha}
	\frac{\zeta(1+i\beta-i\alpha)}{1+i\beta-i\alpha}  - \pi f(\beta) + O(R^{-\delta})\\ 
	&= 2\pi i \sum_{\alpha_B\in\text{LHP}} \mathrm{Res}_{\alpha_B}\left(f(\alpha_B) \left(\frac{T}{2\pi}\right)^{i\beta-i\alpha_B} 
	\frac{\zeta(1+i\beta-i\alpha_B)}{1+i\beta-i\alpha_B}\right)
\end{align*}
where the sum is over the poles $\alpha_B$ of $f$ within the contour in the $\text{LHP}$,  below the pole $\alpha=\beta$ from the zeta function. 
This sum vanishes in the limit $T\to\infty$ thanks to the $T^{i\beta-i\alpha_B}$ factor, since $\Im(\alpha_B) < 0$. Combining the two terms, we get
\begin{equation*}
	\mathrm{PV} \int_{-R}^{R} d\alpha\,  f(\alpha) \left(\zeta(1+i\alpha-i\beta) + \left(\frac{T}{2\pi}\right)^{i\beta-i\alpha}
	\frac{\zeta(1+i\beta-i\alpha)}{1+i\beta-i\alpha}\right)\notag\\
	= 2\pi f(\beta) +  o(1)
\end{equation*}
which gives the desired result in the limit $R, T\to\infty$.
\end{proof} 

\section{Comparison with random matrix theory}\label{sec:CompareRMT}
We now prove Theorem \ref{thm3}, which is the assertion that the kernel $\mathcal{K}_{\text{CUE}}(x,y, N)$ defined in 
Eq.\eqref{eq:CUEKernel} acts as the Dirac delta function on the unit circle, for a suitable class of functions. More precisely, we need to 
establish Eq.\eqref{eq:CUEtheorem} for functions $f$ that are meromorphic in $\mathbb{C}$, holomorphic in the open unit disc $|z| < 1$ and 
satisfy the asymptotic condition $f(z) = o(|z|^{N+1})$ as $|z|\rightarrow \infty$.
 
Using the integral
\begin{equation*}
 	\int_{\text{U}(N)} dM\, \Lambda_{M}(e^{-ix})\, \Lambda_{M^{\dagger}}(e^{iy}) =
 	\frac{1}{1 - e^{i(y-x)}} + \frac{e^{iN(y-x)}}{1 - e^{i(x-y)}}
\end{equation*}
 where $\Lambda_M(z)$ is the characteristic function of the unitary matrices $M \in \text{U}(N)$ (use Eq.(2.16) of Ref.\cite{Conrey_2003}
 for $m=1, n=2, \omega_1=e^{iy}, \omega_2=e^{ix}$),
we write
\begin{equation*}
\int_{0}^{2\pi} dx\, f(e^{ix})\, \mathcal{K}_{\text{CUE}}(x,y,N) = \int_{0}^{2\pi} dx\, f(e^{ix}) \Big(\frac{1}{1 - e^{i(y-x)}} + 
\frac{e^{iN(y-x)}}{1 - e^{i(x-y)}}\Big)
\end{equation*}
For the class of functions as in Theorem \ref{thm3}, consider the integral with $z=e^{ix}$ and $z'=e^{iy}$
\begin{equation*}
\oint_{|z|=1} \frac{dz}{z}\, f(z)\, \left( \frac{z}{z-z'} + \left(\frac{z'}{z}\right)^N \frac{z'}{z'-z} \right)
\end{equation*}
In the first term we integrate along a circular contour approaching the unit disc from inside, i.e., on a circle of radius $|z|=1-\epsilon$ (with 
$\epsilon > 0$ infinitesimally small) so as to avoid the singularity at $z=z'$. Since $f(z)$ is holomorphic inside the unit disc, we obtain
\begin{equation*}
0=\lim_{\epsilon \to 0+} \int_{|z|=1-\epsilon} \frac{dz}{z}\, f(z)\, \frac{z}{z-z'} = i\, \text{PV} \int_{0}^{2\pi} dx\, \frac{f(e^{ix})}{1-e^{i(y-x)}} - i\pi f(z') 
\end{equation*}
as the integrand in the LHS is analytic inside the contour. Similarly, in the second term we integrate along a circular contour approaching 
the unit disc from outside, to get
\begin{equation*}
\lim_{\epsilon \to 0+} \int_{|z|=1+\epsilon} \frac{dz}{z}\, f(z)\, \frac{z'^{N+1} z^{-N}}{z'-z} = i\,\text{PV} \int_{0}^{2\pi} dx\, f(x)\, 
\frac{e^{iN(y-x)}}{1-e^{i(x-y)}} - i\pi f(z')
\end{equation*}
This time, the LHS is not zero. However, we can close the contour around the point at infinity. This picks up contributions from the poles of $f(z)$, 
which may lie outside the unit disc. Let us further assume that $f(z) z^{-N-1}\stackrel{|z|\to\infty}{\longrightarrow} 0$, so that there is no contribution 
from the point at infinity. Note that since all the poles of $f(z)$ are outside the unit disc, the $z^{-N}$ term makes the contribution from the residues
negligibly small as $N \to \infty$. Therefore, the desired result
\begin{equation*}
\int_{0}^{2\pi} dx\, f(e^{ix})\, \mathcal{K}_{\text{CUE}}(x,y,N)  = 2\pi f(e^{iy})
\end{equation*}
is obtained after adding the two contributions. This completes the proof of Theorem \ref{thm3} and establishes the anology with the kernel of the 
second moment of the zeta function with shifted arguments.

We end by noting that the theorem would hold for the polynomials $f_N(e^{ix})=\sum_{0}^N c_k e^{ikx}$, by rewriting the kernel as 
$\mathcal{K}_{\text{CUE}} (x,y,N) = \sum_{0}^{N} e^{in(y-x)}$ and integrating term by term. The proof above is for a more general class 
of functions as in Theorem \ref{thm3}.

\bigskip


\noindent{\bf Acknowledgments:} DG is supported in part by the Mathematical Research Impact Centric Support (MATRICS) grant no.\
MTR/2020/000481 of the Science \&\ Engineering Research Board, Department of Science \&\ Technology, Government of India. 

\bibliographystyle{hieeetr}
\bibliography{RZSSMoment}

\end{document}